\documentclass[twoside,12pt,a4paper]{amsart}
\usepackage{amssymb}
\usepackage{mathrsfs}
\usepackage{amsmath}
 \usepackage[colorlinks,citecolor=blue,urlcolor=blue, linkcolor=blue, backref]{hyperref}
 \usepackage{undertilde}
 \usepackage[capitalise]{cleveref}		
 \usepackage{upgreek}
 \usepackage{tikz}
 
\newtheorem{theorem}{Theorem}[section]
\newtheorem*{theorem*}{Theorem}

\newtheorem{conjecture}[theorem]{Conjecture}
\newtheorem{corollary}[theorem]{Corollary}
\newtheorem{definition}[theorem]{Definition}

\newtheorem{proposition}[theorem]{Proposition}
\newtheorem{question}[theorem]{Question}

\def\uh{\upharpoonright}

\newcommand{\CK}{\omega_1^{\mathrm{CK}}}
\newcommand{\KO}{\mathcal{O}}

\newcommand{\bi}{\begin{itemize}}
\newcommand{\ei}{\end{itemize}}
\newcommand{\bc}{\begin{center}}
\newcommand{\ec}{\end{center}}

\def\uh{\upharpoonright}

\newcommand{\D}{\mathcal{D}}

 \newcommand{\DC}{\mathrm{DC}}

\newcommand{\ZF}{{\mathrm{ZF}}}

\newcommand{\AD}{{\mathrm{AD}}}

\newcommand{\x}{{\mathbf{x}}}

\begin{document}
\title{Some more results on relativized Chaitin's $\Omega$}
 \author{ Liang Yu }
 \thanks{  Yu was partially supported by NSF of China No. 12025103.}
  
\address{Department of Mathematics\ \\
Nanjing University, Jiangsu Province 210093\\
P. R. of China} \email{yuliang.nju@gmail.com.}
 
\maketitle
\begin{abstract}
We prove that, assuming $\ZF$, and restricted to any pointed set, Chaitin's $\Omega_U:x\mapsto \Omega_U^x=\sum_{U^x(\sigma)\downarrow}2^{-|\sigma|}$ is not injective for any universal prefix-free Turing machine $U$, and that $\Omega_U^x$ fails to be degree invariant in a very strong sense, answering several recent questions in descriptive set theory. Moreover, we show that under $\ZF+\AD$, every function $f$ mapping $x$ to $x$-random must be uncountable-to-one over an upper cone of Turing degrees.\end{abstract}
\section{Introduction} \label{sec:introduction}

Given a reduction $\leq_R$ (such as Turing reduction $\leq_T$, or hyperarithmetic reduction $\leq_h$), we say that a perfect set $P$ is {\em $\leq_R$-pointed} if there is a perfect tree $S\subseteq 2^{<\omega}$ so that $[S]=P$ and for any $x\in P$, $S\leq_R x$, where $[S]=\{x\in 2^{\omega}\mid \forall n (x\uh n \in S)\}$. Some times we identify a perfect set with its representation. A perfect set $P$ is {\em uniformly $\leq_R$-pointed} if  there is a perfect tree $S\subseteq 2^{<\omega}$   so that $[S]=P$ and for any $x\in P$, $S$ is $\leq_R$ reducible to $x$ with a fixed reduction. Sometimes we simplify $\leq_T$-pointed as pointed.

Andrew Marks made the following conjecture.
\begin{conjecture}[Marks \cite{PLutz21}]\label{conjecture: marks}
Assume that $\ZF+\DC+\AD$. Given any function $f:2^{\omega}\to 2^{\omega}$, there is a pointed set $[S]$ so that, restricted to $[S]$, either $f$ is constant or $f$ is injective.
\end{conjecture}

Some remarkable conclusions were derived from Conjecture \ref{conjecture: marks}. For example, Patrick Lutz \cite{{PLutz21}} showed that the conjecture implies  that no nonprincipal ultrafilter on the Turing degrees is strictly below Martin measure in the Rudin-Keisler order.

The goal of this paper is to give a counterexample to the conjecture.

We organize the paper as follows: In the first section, we give some notations and basic knowledge; In section 2, we present a general way to refute Conjecture \ref{conjecture: marks} under $\ZF+\AD$; In section 3, we show that $\Omega_U$ is a counterexample under $\ZF$; In the last section, we show that $\Omega_U$ is not degree invariant almost everywhere in the Martin's measure sense.

\bigskip

 We assume that readers have some knowledge of   descriptive set theory,  recursion theory and algorithmic randomness theory.  
\subsection{Set theory}
 
 Our major reference of set theory is \cite{Jech03}.   $\ZF$ is Zermelo-Fraenkel axiom system. $\DC$ is  the axiom of dependent choice. $\AD$ is the axiom of determinacy. Throughout the paper, we work within $\ZF$.
 
A tree $S\subseteq 2^{<\omega}$ is a set downward closed. $[S]$ is the collection of infinite paths trough $S$. Given any $x\in \omega^{\omega}$ and natural number $n$, we use $x\uh n$ to denote an initial segment of $x$ with length $n$.  In other words, $x\uh n$ is a finite string $\sigma\in \omega^{<\omega}$ of length $n$ so that for any $i<n$, $\sigma(i)=x(i)$. 

$\langle \cdot, \cdot \rangle$ is G\" odel's pairing function.

The following fundamental  result is due to Martin.
\begin{theorem}[Martin \cite{Martin68}]\label{theorem: std theorem}
If $A\subseteq 2^{\omega}$ so that $\forall x\exists y\geq_T x y\in A$, then $A$ has a pointed subset.
\end{theorem}
 
\subsection{Recursion theory} 

 Our major references of recursion theory are \cite{Sacks90} and \cite{Ler83}.
We use $\leq_T$ to denote Turing reduction and $\leq_h$ to denote hyperarithmetic reduction.  We use $\Phi^x$ denote a Turing machine with oracle $x$. Sometimes we also say that $\Phi^x$ is a recursive functional. We fix an effective enumeration $\{\Phi_e^x\}_{e\in \omega}$ of recursive functionals.

Given a real $x$, its    Turing degree $\x$ is a set of reals defined as $\{y\mid y\equiv_T x\}$. We say $\x\leq \mathbf{y}$ if $x\leq_T y$.  We use $\D$ to denote the set of Turing degrees. An {\em upper cone}   of Turing degrees is a set $\{\mathbf{y}\mid \mathbf{y}\geq \mathbf{x}  \}$ for some fixed $\mathbf{x}$.

\begin{theorem}[Jockusch and Simpson \cite{JoSi76}]\label{theorem: uniformly subset}
Any $\leq_T$ pointed set $P$ has a uniformly $\leq_T$-pointed subset.
\end{theorem} 

 Kleene's $\KO$, which is a standard $\Pi^1_1$-complete set, is as defined in \cite{Sacks90}. $\CK$ is the least non-recursive ordinal and $\omega_1^x$ is the least ordinal not recursive in $x$.

We say that a set $A$ ranges Turing degrees cofinally if for any real $x$, there is some $y\geq_T x$ in $A$. We use $x'$ to denote the Turing jump relative to $x$. More generally, if $\alpha<\omega_1^x$, then $x^{(\alpha)}$ is that $\alpha$-th Turing jump of $x$.

\begin{theorem}[Martin \cite{Martin76}]\label{theorem: friedman's conjecture}
For any $\Delta^1_1$  set $A\subseteq 2^{\omega}$, if $A$ has a nonhyperarithmetic real, then it  has a $\leq_h$-pointed subset.
\end{theorem}

\subsection{Algorithmic randomness theory}

For the classical algorithmic randomness theory, see  \cite{Niesbook09} and \cite{DH10}. A prefix-free Turing machine $M$ is a Turing machine so that for any  input $\sigma\in 2^{<\omega}$, if $M(\sigma)$ halts, then $M(\tau)$ does not for any $\tau$ extending $\sigma$. A universal prefix-free Turing machine is a prefix-free Turing machine coding all the others.

A real $x$ is left-r.e. if there is a recursive increasing sequence rationals $\{q_s\}_{s\in \omega}$ so that $\lim_{s\to \infty} q_s=x$.

\begin{theorem}[Chaitin \cite{Cha75}]\label{theorem: chaitin theorem}
Given any universal prefix-Turing machine $U$, $\Omega_U^x=\sum_{U(\sigma)\downarrow}2^{-|\sigma|}$ is a left-r.e. Martin-L\" of random real relative to $x$.  
\end{theorem}

Sometimes, we simply say random instead of Martin-L\" of random, and $x$-random instead of Martin-L\" of random relative to $x$.

\begin{theorem}[van Lambalgen \cite{van90}]\label{theorem: van lambalgen}
For any random reals $x$ and $y$, $x$ is $y$-random if and only if $y$ is $x$-random if and only if $x\oplus y$ is random.
\end{theorem}

\begin{theorem}[Ku\v{c}era \cite{Ku86}]\label{theorem: kucera}
If $x$ is random and left-r.e., then $x\equiv_T \emptyset'$, the Turing jump of $\emptyset$.
\end{theorem}

A notion so-called {\em Schonrr random} is randomness notion weaker than \break Martin-L\" of random.  A real $x$ is called {\em low for Schonorr-random} if every Schnorr random is Schnorr random relative to $x$.

\begin{theorem}[Terwijn and Zambella \cite{TZ}]\label{theorem: zt}
There is a perfect tree $S\leq_T \emptyset''$ only containing low for Schnorr random reals.
\end{theorem}

For the facts of    higher randomness theory,  see \cite{CY15book}.  A real $r$ is $\Pi^1_1$-random if it does not belong to any $\Pi^1_1$-null set.

\section{ Solovay's uniformization principle and its application}
The following uniformization principle is due to Solovay (unpublished).
\begin{theorem}[Solovay]\label{theorem: solovay uniformization}
Assuming $\ZF+\AD$, for any binary relation $R\subseteq 2^{\omega}\times 2^{\omega}$ with $\forall x\exists y R(x,y)$, there is a Borel function $g:2^{\omega}\to 2^{\omega}$ so that  the set $\{x\mid R(x,g(x))\}$ is conull.
\end{theorem}
We have not found a proof of the theorem. The following proof is  recursion theoretical.
\begin{proof}
 Let $$\tilde{R}=\{(x,y)\mid \forall z_0\leq_T x\exists z_1\leq_T y R(x,y)\}.$$  By Theorem \ref{theorem: zt}, for every $x$, there is some $y>_T x$ low for $x$-schnorr random so that $y''$ computes an $x$-Schnorr random $r$ and a real $z$ so that  $\tilde{R}(r,z)$. So there are numbers $e_0$ and $e_1$ so that the set $$A=\{y\mid \tilde{R}(\Phi_{e_0}^{y''}, \Phi_{e_1}^{y''})\wedge \Phi_{e_0}^{y''}\mbox{ is Schnorr random relative to }y\}$$ ranges Turing degrees cofinally. By Theorem \ref{theorem: std theorem}, $A$ has a pointed subset $S$. Let $$B=\{r\mid \exists x\in [S] (r=\Phi^{x''}_{e_0})\}.$$ Then $B$ has a positive measure.   Define $$\hat{R}=\{(r,z)\mid \exists y\in S(r=\Phi_{e_0}^{y''} \wedge z=\Phi_{e_1}^{y''}) \}.$$ Then $\hat{R}\subseteq \tilde{R}$ is $\Sigma^1_1(S)$. Then there is a Borel function $f$ and positive measure Borel set $C\subseteq B$ so that for any $x\in C$, $\hat{R}(x,f(x))$ and so $\tilde{R}(x,f(x))$.   Let $$[C]_T=\{x\mid \exists y\in C(x\equiv_T y)\}.$$ Then $[C]_T$ is a Borel conull set. For any $r\in [C]_T$,  let $\langle e_0,e_1\rangle$ be the least pair so that $\Phi_{e_0}^r\in C$ and $\Phi_{e_1}^{\Phi_{e_0}^r}=r$. Then $f(\Phi_{e_0}^r)$ is defined and $\tilde R(\Phi_{e_0}^r,  f(\Phi_{e_0}^r))$. By the definition of $\tilde{R}$ and the fact that $r\leq_T \Phi_{e_0}^r$, we may let $e_2$ be the least index so that $R(r, \Phi_{e_2}^{  f(\Phi_{e_0}^r)})$.
 
 So for any pair $\langle e_0,e_2\rangle$,  define 
 
 \bigskip
 $g_{\langle e_0,e_2\rangle}(r)= \left\{
\begin{array}{r@{\quad\quad}l}
\Phi_{e_2}^{  f(\Phi_{e_0}^r)}, & \Phi_{e_2}^{  f(\Phi_{e_0}^r)} \mbox{ is defined};  \\
r, & \mbox{Otherwise.}\\
\end{array}
\right. $

 Set $$E_{\langle e_0,e_2\rangle }=\{r\mid \Phi_{e_2}^{  f(\Phi_{e_0}^r)}\mbox{ is defined }\wedge R(r, \Phi_{e_2}^{  f(\Phi_{e_0}^r)})\}.$$Then $g_{\langle e_0,e_2\rangle}$ is a Borel function. By the discussion above,  we have that $[C]_T\subseteq \bigcup_e E_{e }$ and so $\bigcup_{e}E_{e }$ is conull. Since $\AD$ implies that every set is measurable, we have that for any $e$, there is a Borel set $F_{e}\subset E_{e }\setminus \bigcup_{   e'  < e }E_{e' }$ so that    $E_{e }\setminus( F_{e}\cup \bigcup_{ e'<e }E_{e' })$ is null. Clearly $\{F_e\}_{e\in \omega}$ is a disjoint family so that $\bigcup_{e}F_{e }$ is conull.
\bigskip
 
 Define
 $g(x)= \left\{
\begin{array}{r@{\quad\quad}l}
\Phi_{e_2}^{  f(\Phi_{e_0}^x)}, & \exists e_0\exists e_2(x\in F_{\langle e_0,e_2\rangle});  \\
x, & \mbox{Otherwise.}\\
\end{array}
\right. $

\bigskip

 Then $g$ is a well-defined Borel function and  for almost every real $x$, $R(x,g(x))$.
\end{proof}

\begin{corollary}\label{corollary: boundedness}
If $R\subseteq 2^{\omega}\times \omega_1$ so that $\forall x\exists \alpha R(x,\alpha)$, then there is  some $\alpha_1<\omega_1$ so that the set $\{x\mid \exists \alpha\leq \alpha_1 R(x,\alpha)\}$ is conull.
\end{corollary}
\begin{proof}
Define $\tilde R\subseteq 2^{\omega}\times 2^{\omega}$ so that $\tilde R=\{(x,z)\mid \exists \alpha\leq \omega_1^z R(x,\alpha)\}$. 

By Theorem \ref{theorem: solovay uniformization}, there is a Borel function $f: 2^{\omega}\to 2^{\omega}$ so that the set \break $\{x\mid \tilde R(x,f(x))\}$ is conull. Let $F$ be a Borel conull set so that $\forall x\in F \tilde R(x,f(x))$. Then there is a real $z$ so that both $f$ and $F$ are hyperarithmetic in $z$.  Then for any $\Pi^1_1(z)$-random real $r\in F$, $f(r)\leq_h r\oplus z$ and so $\omega_1^{f(r)}\leq \omega_1^{r\oplus z}=\omega_1^z$. Let $\alpha_1=\omega_1^z$.
\end{proof}

\begin{definition}
For any function $f:2^{\omega}\to 2^{\omega}$ and $y\in 2^{\omega}$, define $$\alpha^f_y=\sup\{\omega_1^x\mid f(x)=y\}.$$
\end{definition}

 



The following theorem provide a general way to refute Conjecture \ref{conjecture: marks} under $\ZF+\AD$. 

\begin{theorem}\label{theorem: general theorem}
Assume $\ZF+\AD$. If $f:2^{\omega}\to 2^{\omega}$ is a function so that for any $x$, $f(x)$ is Martin-L\" of random relative to $x$, then there is real $x_0$ so that for any $x\geq_T x_0$, $\alpha^f_{f(x)}=\omega_1$.
\end{theorem}
\begin{proof}
 Suppose not, then   the set   $\{x   \mid \alpha^f_{f(x)}<\omega_1\}$ ranges Turing degrees cofinally. Then by Martin's theorem \ref{theorem: std theorem}, $\{x \mid \alpha^f_{f(x)}<\omega_1\}$ has a pointed subset $S$. 
 
 Clearly the range of $f$ over $S$ has positive measure.
 
 Now for any countable ordinal $\alpha$, let $$A_{\alpha}=\{x\in S\mid \alpha^f_{f(x)}=\alpha\} \mbox{ and }B_{\alpha}=\{f(x)\mid x\in A_{\alpha}\}.$$ Then $\{B_{\alpha}\}_{\alpha<\omega_1}$ is a disjoint family so that $\forall \alpha \mu(\bigcup_{\beta>\alpha}B_{\beta})>0$. 
 
 Define a binary relation $R\subseteq 2^{\omega}\times \omega_1$ so that $R=\{(y,\alpha)\mid y\in B_{\alpha}\}$. By Corollary \ref{corollary: boundedness}, there is some $\alpha_1<\omega_1$ so that the set $\{y\mid \exists \alpha\leq \alpha_1 R(y,\alpha)\}$ is conull. Then the set $\bigcup_{\alpha>\alpha_1}B_{\alpha}$ is null. This is a contradiction. 
 \end{proof}

An immediate conclusion of Theorem \ref{theorem: general theorem} is that $\Omega_U$ operator is a counterexample to the Conjecture \ref{conjecture: marks}. 

\begin{corollary}Assume $\ZF+\AD$.
Given any universal prefix-Turing machine $U$, there is a real $x_0$ so that for any $x\geq_T x_0$, $\{y \mid \Omega_U^y=\Omega_U^x\}$ has a $\leq_h$-pointed subset. So $\Omega^x_U$ is a counterexample to Conjecture. \footnote{Lutz proves that a corrected version of  Conjecture \ref{conjecture: marks}. I.e.  for any $f$, if for any pointed set $S$, $f(S)$ ranges Turing degrees cofinally, then $f$ is injective on some pointed set.} 
\end{corollary} 
\begin{proof}
Clearly $\Omega_U$ satisfies the assumption of Theorem \ref{theorem: general theorem}. So there is a real $x_0$ so that for any $x\geq_T x_0$, the set  $\{y\in S \mid \Omega_U^y=\Omega_U^x\}$ is uncountable. Note that   $\{y\in S \mid \Omega_U^y=\Omega_U^x\}$  is  $\Delta^1_1(\Omega_U^x)$ and only contains real hyperarithmetically above $\Omega_U^x$. By Theorem \ref{theorem: friedman's conjecture}, it has a  $\leq_h$-pointed subset.
\end{proof}
 
\section{The counterexample under $\ZF$}

We  present  two methods to show, under $\ZF$, that $\Omega^x_U$ is a  counterexample to Conjecture \ref{conjecture: marks}. 

The first method is an application of  higher randomness theory.   

\begin{theorem}\label{theorem: main 1}
For any $\leq_T$-pointed set $S$ and universal prefix-free Turing machine $U$,  there is a real $z\in  S $ so that $\{y\in S \mid \Omega_U^y=\Omega_U^z\}$ has a $\leq_h$-pointed subset. So $\Omega^x_U$ is a counterexample to Conjecture \ref{conjecture: marks}.
\end{theorem}
\begin{proof}
We first prove that $\Omega^x_U$ cannot be constant over a $\leq_T$-pointed set. Suppose not. Then $\Omega^x_U=r$ for some real $r$  over a $\leq_T$-pointed set.  Then there is some real $z_0$ such that $z_0\geq_T r$ and $\Omega^{z_0}_U=r$. But $r$ is a $z_0$-random and so is not Turing reducible to $z_0$, which is a contradiction.

  Now fix a real $z\equiv_T \KO^{(\KO^S)}$ in $S$. Such a real $z$ exists since $S$ is a $\leq_T$-pointed set. Then $\Omega_U^z$ is $\Delta^1_1(\KO^S)$- and so $\Pi^1_1(S)$-random (see \cite{CNY08}). Also note that $$\Omega_U^z\oplus S\leq_h z.$$ Let $$A=\{y\in S\mid \Omega_U^y=\Omega_U^z\}.$$ It is clear that $A$ is $\Delta^1_1(\Omega_U^z\oplus S)$ and $z\in A$. Since $\Omega_U^z$ is $\Pi^1_1(S)$ random, we have that $\omega_1^{\Omega_U^z\oplus S}=\omega_1^S$ by Sacks \cite{Sacks69} (also see Corollary 14.3.2 in \cite{CY15book}).  Thus we have that  $ z\not\leq_h \Omega^z\oplus S$ and so $$\Omega_U^z\oplus S<_h z.$$ Also since $S$ is pointed, we have that $$\Omega^z_U\oplus S=\Omega^y_U\oplus S\leq_h y$$ for any $y\in A$. Then by Theorem \ref{theorem: friedman's conjecture} relative to  $\Omega_U^z\oplus S$, we have that $A$ has a $\leq_h$-pointed subset.
\end{proof}
 
A natural question is (by Lutz \cite{PLutzemail}) that given any $\leq_T$-pointed set $S$, are there two different reals $x\equiv_T y$ in $S$ so that $\Omega_U^x=\Omega_U^y$? We give a positive answer to the question via classical randomness theory and so the second way to refute Conjecture \ref{conjecture: marks}. The basic ideas are from \cite{DHMN} and \cite{RS:08}. 

Fix a pointed set $S$ and universal prefix-free Turing machine $U$.  Clearly there is an $S$-recursive  function $f:2^{<\omega} \to 2^{<\omega}$ so that $\hat{f}: x\mapsto \bigcup_n f(x\uh n)$ is a homeomorphism from $2^{\omega}$ to $S$. Let $\hat{S}_1\subset 2^{<\omega}$ be an $S$-recursive tree so that
\begin{itemize}
\item  $\mu(S_1)>0$; and 
\item  $S_1$ only contains $S$-random reals. 
\end{itemize}

Set $$r=\inf\{\Omega_U^{\hat{f}(x)}\mid x\in S_1\}.$$ Clearly $r$ is left-r.e. relative to  $S$ and so $r\oplus S\equiv_T (S)'$. Since $S_1$ is compact, there is a real $x_0\in S_1$ so that $r=\Omega_U^{\hat{f}(x_0)}$.\footnote{This can be proved by a standard technique in recursion theory but may not be   so clear to the readers not familiar with it. The proof of Theorem 6.1 in \cite{DHMN}   contains more details.  } Since $S$ is pointed, we have that $$\hat{f}(x_0)\equiv_T S\oplus x_0.$$ So $r$ is random   relative to $S\oplus x_0$ and so to $S$. By van Lambalgen Theorem \ref{theorem: van lambalgen} and Theorem \ref{theorem: kucera} relative to $S$, we have that $x_0$ is random relative to $S\oplus r\equiv_T (S)'$. Since $r$ is left-r.e. relative to  $S$,    the set $$F^U_r=\{y \mid \Omega_U^{\hat{f}(y)}=r\}$$ is $\Pi^0_2(S)$ and so measurable. Since $x_0\in F^U_r$ is $S'$-random,  $F^U_r$ must have positive measure.

So there must be some reals $y_0,y_1\in F^U_r$ so that $y_0\equiv_T y_1$ (otherwise,  $F^U_r$ must be null).    Then $\hat{f}(y_0)\oplus S\equiv_T y_0\oplus S\equiv_T y_1 \oplus S\equiv_T \hat{f}(y_1)\oplus S$ since $f\leq_T S$  and $\hat{f}$ is a homeomorphism. By the pointedness of $S$, we have that $\hat{f}(y_0)\geq_T S$ and $\hat{f}(y_1)\geq_T S$. Thus  $$\hat{f}(y_0) \equiv_Ty_0\oplus S\equiv_T y_1 \oplus S\equiv_T \hat{f}(y_1).$$  

So $$f(y_0)\equiv_T f(y_1),\Omega^{f(y_0)}_U=\Omega^{f(y_1)}_U \mbox{ but } f(y_0)\neq f(y_1).$$

Hence we have the following theorem.
\begin{theorem}\label{theorem: classcial proof}
Given any universal prefix-free Turing machine $U$ and any pointed set $S$, there are two different reals $z_0$ and $z_1$ in $S$ so that $\Omega^{z_0}_U=\Omega^{z_1}_U$ and $z_0\equiv_T z_1$.
\end{theorem}

\section{On degree invariantness}

A function $f:2^{\omega}\to 2^{\omega}$ is {\em degree invariant} if $\forall x\forall y(x\equiv_T y \implies f(x)\equiv_T f(y))$. The following question is open to us.
\begin{question}\label{question: degree invaraint}
Is it consistent with $\ZF+\DC$ that there is a  degree invariant Borel function $f$ so that for any $x$, $f(x)$ is random relative to $x$?
\end{question}

Clearly Question \ref{question: degree invaraint} is related to Martin's conjecture.\footnote{\cite{PLutz21} contains an up-to-date survey concerning Martin's conjecture.} One may wonder whether Chaitin's $\Omega$ can be served as a solution to the question. But in \cite{DHMN}, it has been shown that $\Omega_U^x$ is not degree invariant for any universal prefix-free Turing machine $U$. The following result says that $\Omega_U^x$ fails to be degree invariant in a very strong sense.
\begin{proposition}\label{proposition: turing incomparable}
Given any universal prefix-free Turing machine $U$ and pointed set $S$, there are two reals $z_0,z_1\in S$ so that $z_0\equiv_T z_1$ but $\Omega^{z_0}_U$ is Turing incomparable with $\Omega^{z_1}_U$.
\end{proposition}
\begin{proof}

The proof is based on the proof of Theorem 6.7 in \cite{DHMN} via a pushout-pullback method.  We follow the notations in the proof of Theorem \ref{theorem: classcial proof}.

By Theorem \ref{theorem: uniformly subset}, we may assume that $S$ is uniformly pointed. So we may assume that there is a recursive function $\Phi$ so that for any $x\in S$, $\Phi^x=S$. Then there is another recursive functional $\Psi$ so that for any $x\in S$, $\Psi^x=f^{-1}(x)$. Now let $V$ be another universal prefix-free Turing machine so that for any $\sigma\in 2^{<\omega}$ and real $x$, $V^x(0\sigma)=U^x(\sigma)$ and $V^x(1\sigma)=M^{\Psi^x}(\sigma)$, where $M^{z}$ is a prefix-free Turing machine so that $\Omega_{M }^z=\sum_{M^z(\sigma)\downarrow}2^{-|\sigma|}=z$. Then for any real $x\in S$, $$\Omega_{V}^x=\frac{\Omega_U^x+\hat{f}^{-1}(x)}{2}.$$

  By replacing $U$ with $V$ in the proof of Theorem \ref{theorem: classcial proof}, there is another left-r.e real $r_0$ relative to $S$ so that $$F^V_{r_0}=\{y\in 2^{\omega}\mid \Omega_V^{\hat{f}(y)}=r_0\}$$ has positive measure. Let $z\geq_T S'$ and $y \in F^V_{r_0}$ be any real random relative to $z$, then $$\Omega_U^{\hat{f}(y )}=2 \Omega_V^{\hat{f}(y)}-\hat{f}^{-1}(\hat{f}(y ))=2r_0-y .$$

Since $r_0\leq_T S'\leq_T z$ and $y$ is $z$-random, we have that $\Omega_U^{\hat{f}(y )}$ must also be $z$-random. So the set $$G_z^U=\{y\mid \Omega_U^{\hat{f}(y)}\mbox{ is }z\mbox{-random}\}$$ has positive measure.

Now let $$\hat{F}^U_r=\{y^*\mid \exists y\in F^U_r (y\equiv_T y^*)\}\mbox{ and } \hat{G}^U_z=\{y^*\mid \exists y\in G^U_z(y\equiv_T y^*)\}.$$

Then both the sets have measure 1 and so there must be   $s_0\in \hat{F}^U_r$ and $s_1\in \hat{G}^U_r$ so that $s_0\equiv_T s_1$. Then we have that $$\hat{f}(s_0) \equiv_Ts_0\oplus S\equiv_T s_1 \oplus S\equiv_T \hat{f}(s_1).$$ Since  $z\geq_T S'\geq_T r$,  we have that  $\Omega_U^{\hat{f}(s_1)}$ must be $r$-random and so  Turing incomparable with $r=\Omega_U^{\hat{f}(s_0)}$ by van Lambalgen Theorem \ref{theorem: van lambalgen}. Set $z_0=f(s_0)$ and $z_1=f(s_1)$. They are as required. 
\end{proof}

By combining Theorem \ref{theorem: classcial proof},  Proposition \ref{proposition: turing incomparable}  and Borel determinacy \cite{Martin75}, we have the following conclusion.
\begin{corollary}
Assuming $\ZF+\DC$ \footnote{$\DC$ is the axiom of dependent choice.},  there is a real $x$ so that for any universal prefix-free Turing machine $U$ and any real $y\geq_T x$, there are three different reals $z_0\equiv_T z_1 \equiv_T z_2\equiv_T y$ so that  $\Omega_U^{z_0}= \Omega_U^{z_1}$ is Turing incomparable with $\Omega_U^{z_2}$.
\end{corollary}
\begin{proof}
We follow the notations in the proof of Proposition \ref{proposition: turing incomparable}. Given any pointed set $S$ and universal prefix-free Turing machine $U$, by the proofs of Theorem  \ref{theorem: classcial proof} and Proposition \ref{proposition: turing incomparable} (we assume $z\geq_T S'\geq_T r$), both $F^U_r$ and $G^U_z$ have positive measure. Then it is clear that the Borel set $$\{y\in F^U_r  \mid \exists y^*\equiv_T y(y^*\neq y\wedge y^*\in F^U_r)\}$$ must also have positive measure. So there must be three different reals $y_0\equiv_T y_1\equiv_T y_2$ so that $y_0,y_1\in  F^U_r$ and  $y_2\in  G^U_z$. So $f(y_0)\equiv_T f(y_1)\equiv_T f(y_2)\geq_T S$ are three different reals so that $\Omega_U^{f(y_0)}= \Omega_U^{f(y_1)}=r$ is Turing incomparable with $\Omega_U^{f(y_2)}$. Thus the set 
   \begin{multline*}B_U=\{y\mid \mbox{There are three different reals }z_0\equiv_T z_1 \equiv_T z_2\equiv_T y\mbox{ so that  }\\ \Omega_U^{z_0}= \Omega_U^{z_1}\mbox{ is Turing incomparable with }\Omega_U^{z_2}\}\end{multline*} is a Borel set of Turing degrees that ranges Turing degrees cofinally. By the Borel determinacy, $B_U$ contains an upper cone of Turing degrees. But there are only countably many such Turing machines. So $$\bigcap\{B_U\mid U\mbox{ is a universal prefix-free Turing machine}\}$$ contains an upper cone of Turing degrees.
\end{proof}

\bigskip

\bibliographystyle{plain}

\end{document}